\documentclass[12pt, reqno]{amsart}
\usepackage{ifthen,amsfonts,
amsmath,amssymb, enumerate,
epsfig,color,graphicx,bbm,graphics}
\usepackage{stmaryrd}
\usepackage{caption}
\usepackage{a4wide}
\usepackage{xfrac}
\usepackage{mathrsfs,bm,amsthm,mathtools,yfonts}
\usepackage{xcolor}
\usepackage{tikz}
\usepackage{float}
\usepackage{subfig}

\mathtoolsset{showonlyrefs,showmanualtags}

\makeatletter

\begingroup

\newtheorem{theorem}{Theorem}[section]

\endgroup

\makeindex
\begin{document}

\title[Smoothing and selection for transport equations]
{Smoothing does not give a selection principle for transport equations with bounded autonomous fields}
\author{Camillo De Lellis}
\address{School of Mathematics, Institute for Advanced Study, 1 Einstein Dr., Princeton NJ 05840, USA\\
and Universit\"at Z\"urich}
\email{camillo.delellis@math.ias.edu}

\author{Vikram Giri}
\address{Department of mathematics, Princeton Uniiversity, Washington Rd., Princeton NJ 05840, USA}
\email{vgiri@math.princeton.edu}

\begin{abstract}
We give an example of a bounded divergence free autonomous vector field in $\mathbb R^3$ (and of a nonautonomous bounded divergence free vector field in $\mathbb R^2$) and of a bounded initial data for which the Cauchy problem for the corresponding transport equation has $2$ distinct solutions. We then show that both solutions are limits of classical solutions of transport equations for appropriate smoothings of the vector fields and of the initial data.  
\end{abstract}

\maketitle


\section{Introduction}
In this note, we consider the classical Cauhy problem for a transport equation of type
\begin{equation}\label{e:TE}
\begin{cases}
    \partial_t \theta + (v \cdot \nabla_x) \theta = 0 \\
    \theta (0, x) = \theta_{in}(x)\, 
\end{cases}
\end{equation}
on $[0,T]\times \mathbb R^{d+1}$ (with $d\geq 2$), where $\theta$ is the unknown, while $v$ is a known vector field.
The vector fields considered will be divergence free and thus \eqref{e:TE} can be rewritten as $\partial_t \theta + {\rm div}\, (v\theta) =0$ (which is usually called continuity equation). Therefore, as it is customary in  the literature, 
when $v$ and $\theta$ are summable enough (i.e. $v\in L^p$ and $\theta\in L^{p'}$ for a pair of dual exponents $p, p'$) we understand solutions in the distributional sense. 

We will restrict our attention to initial data $\theta_{in}$ which are bounded, to solutions which are bounded and to vector fields which are bounded. Under such assumptions \eqref{e:TE} is classically well-posed if $v$ is a Lipschitz vector field. Moreover the solutions are stable for perturbations of the vector field $v$. The famous seminal paper \cite{DL} established a similar well-posedness and stability theory when $v\in L^1 ([0,T], W^{1,p} (\mathbb R^n))$ for any $p\in [1, \infty]$: this is commonly called DiPerna-Lions theory and it has far-reaching applications to very different problems. The DiPerna-Lions theory was extended by Ambrosio in \cite{A04} to $L^1 ([0,T], BV (\mathbb R^n))$ and it was then showed that the result is essentially optimal: weak solutions for vector fields $v\in W^{s,1}$ are in general not unique for $s<1$ (cf. \cite{Aizenman,D}; nonetheless there are still several important open problems in the area and very recent interesting developments, see for instance \cite{ACoFi,AlBiCr,AlCrMa,BiBo,BouCr,BrCoDe,ChJa,CoTi,CrDe,MoSz,YaZl}). 

The next natural question in this regard is then whether there is a meaningful selection principle among these different weak solutions. For instance, do solutions of suitable regularizations have a unique limit? To our knowledge this question is specifically raised for the first time in \cite{CCS20}, where the authors give a partial negative answer. 
The aim of this note is to show that, at least if we only require the regularizations to just enjoy (in a uniform way) the same regularity estimates of the vector field, then the answer is negative. The answer is negative even if we consider autonomous vector fields and if the initial data remains fixed, or anyway they are regularized by convolution with a classical kernel.
Our main theorem is the following:

\begin{theorem}\label{t:main} Let $d\geq 2$. Then
    there exist
    \begin{itemize}
    \item[(i)] an autonomous compactly supported divergence-free vector field $v \in L^{\infty}(\mathbb{R}^{d+1} ; \mathbb{R}^{d+1})$,
    \item[(ii)] an initial data $\theta_{in} \in L^{\infty}(\mathbb{R}^{d+1};\mathbb{R})$ with compact support,
    \item[(iii)] two sequences of divergence-free vector fields $\{ v_i' \}_{i=1}^{\infty}, \{ \tilde v_i \}_{i=1}^{\infty} \subset C^{\infty}_c ( \mathbb{R}^{d+1} ; \mathbb{R}^{d+1})$, 
    \item[(iv)] and a sequence of smooth initial data $\{\theta_{i, in}\}_{i=1}^{\infty} \subset C^\infty_c(\mathbb{R}^{d+1}; \mathbb{R})$, 
\end{itemize}
 with the following properties.  If $\theta_i'$ and $\tilde \theta_i$ are the unique solutions to the transport equation (TE) with initial data $\theta_{i, in}$, then:
    \begin{itemize}
        \item[(a)] $v_i' \to v$, $\tilde v_i \to v$, and $\theta_{i,in}\to \theta_i$ strongly in $L^1$ as $i \to \infty$, 
        \item[(b)] $\|v_i'\|_{L^\infty}, \|\tilde v_i\|_{L^\infty}, \|\theta_{i, in}\|_{L^\infty} \leq C$ for some constant $C$ independent of $i$,
        \item[(d)] $\theta'_i \rightharpoonup \theta'$ and $\tilde \theta_i \rightharpoonup \tilde \theta$ weakly in $L^1$, where $\theta'$ and $\tilde \theta$ are 2 \underline{distinct} solutions to \eqref{e:TE}.
        \end{itemize}
Moreover, the vector field $v$ belongs to $W^{s,p}$ for every $s<1$ and $p<\frac{1}{s}$ and the regularized fields $v_i', \tilde{v}_i$ enjoy uniform estimates in the corresponding spaces, while the initial data is piecewise constant, has bounded variation, and is regularized to $\theta_{i,in}$ through convolution with a standard kernel.
\end{theorem}

Previous work of \cite{CCS20} has shown this theorem for a suitable field $v \in L^p$ for $p \in [0, \frac{4}{3}]$, with a completely different construction. 
If we drop the requirement that the field be autonomous, we can show the same conclusion for $2$-dimensional fields.

\begin{theorem}\label{t:main-2} There exist
    \begin{itemize}
    \item[(i)] a compactly supported divergence-free vector field $b \in L^{\infty}([0,T]\times \mathbb{R}^{2} ; \mathbb{R}^{2})$,
    \item[(ii)] an initial data $\rho_{in} \in L^{\infty}(\mathbb{R}^{2};\mathbb{R})$ with compact support,
    \item[(iii)] two sequences of divergence-free vector fields $\{ b_i' \}_{i=1}^{\infty}, \{ \tilde b_i \}_{i=1}^{\infty} \subset C^{\infty}_c ([0,T]\times \mathbb{R}^{2} ; \mathbb{R}^{2})$, 
    \item[(iv)] and a sequence of smooth initial data $\{\rho_{i, in}\}_{i=1}^{\infty} \subset C^\infty_c(\mathbb{R}^{2}; \mathbb{R})$, 
 with the following properties.
 \end{itemize}  
 If $\rho_i'$ and $\tilde \rho_i$ are the unique solutions to the transport equation \eqref{e:TE} with fields $b'_i$ and $\tilde{b}_i$ and initial data $\rho_{i, in}$, then:
    \begin{itemize}
        \item[(a)] $b_i' \to b$, $\tilde b_i \to b$, and $\rho_{i, in} \to \rho_{in}$ strongly in $L^1 ([0,T]\times \mathbb R^2)$ as $i \to \infty$, 
        \item[(b)] $\|b_i'\|_{L^\infty}, \|\tilde b_i\|_{L^\infty}, \|\rho_{i, in}\|_{L^\infty} \leq C$,
        \item[(c)] $\rho'_i \rightharpoonup \rho'$ and $\tilde \rho_i \rightharpoonup \tilde \rho$ weakly in $L^1$, where $\rho'$ and $\tilde \rho$ are 2 \underline{distinct} solutions to \eqref{e:TE} with field $b$ and initial data $\rho_{in}$.
        \end{itemize}
\end{theorem}

Indeed, since there is a simple way to pass from a non-autonomous example to an autonomous one in one dimension higher, we will mainly focus on how to build the example of Theorem \ref{t:main-2}. The construction is similar to other ones present in the literature, starting from the work of DePauw \cite{D}: the contribution of this note is to show how it can be arranged so that the corresponding distinct solutions are limits of solutions of appropriate regularizations. 

Even though a ``closure'' of classical solutions does not provide a selection mechanism to single out one preferred solution to the final transport equation, our construction does not rule out the possibility that some ``canonical'' regularization (like smoothing by convolution with some specific kernel) still selects only one preferred solution in the limit. 

\subsection{Acknowledgments} C.D.L. acknowledges the support of the NSF grants DMS-1946175 and DMS-1854147, while V.G. acknowledges the support of the NSF grant DMS-FRG-1854344.

\section{Construction of the nonautonomous field} 

In this section we detail the construction of the vector fields and of the initial data in Theorem \ref{t:main}. We only ignore two aspects: the vector fields (and the initial data) will not be compactly supported and we do not give estimates on their $W^{s,p}$ norms. Both aspects are minor and will be addressed in the Section \ref{s:minor}, where we also show how to pass to the autonomous example in one dimension higher. 

\subsection{Step 1.  Definition of $b$, $\rho_{in}$, $\rho'$, and $\tilde{\rho}$.}
We first introduce the following two standard lattices on $\mathbb R^2$, namely $\mathcal{L}^1 := \mathbb Z^2\subset \mathbb R^2$ and $\mathcal{L}^2:=\mathbb Z^2 + (\frac{1}{2}, \frac{1}{2})\subset \mathbb R^2$. To both of them we can associate a corresponding subdivision of the plane into squares which have vertices lying in the corresponding lattices, which we denote by $\mathcal{S}^1$ and $\mathcal{S}^2$. We then consider the rescaled lattices $\mathcal{L}^1_i:= 2^{-i} \mathbb{Z}^2$ and $\mathcal{L}^2_i := (2^{-i-1},2^{-i-1})+2^{-i} \mathbb Z^2$ and the corresponding square subdivision of $\mathbb Z^2$, respectively $\mathcal{S}^1_i$ and $\mathcal{S}^2_2$. Observe that
\begin{itemize}
\item[(D)] The centers of the squares $\mathcal{S}^1_i$ are elements of $\mathcal{L}^2_i$ and viceversa.
\end{itemize}
We let $\rho_{in} (x) = \lfloor{x_1}\rfloor + \lfloor{x_2}\rfloor \ mod \ 2$. This is a `chessboard' pattern based on the standard lattice $\mathbb{Z}^2 \subset \mathbb{R}^2$: if we index the squares of $\mathcal{S}^1$ with $(k,j)$, where $(k+\frac{1}{2}, j+\frac{1}{2})\in \mathcal{L}^2$ is the center of the corresponding square, then $\rho_{in}$ vanishes on the squares for which $k+j$ is even, while it is identically equal to $1$ on squares for which $k+j$ is odd.
 
 Next we define the following $2$-dimensional vector field:
 \[
 w(x) =
    \begin{cases}
        (0, 4x_1)\text{ , if }1/2 > |x_1| > |x_2| \\
        (-4x_2, 0)\text{ , if }1/2 > |x_2| > |x_1| \\
        (0, 0)\text{ , otherwise.} \\
    \end{cases}
    \]
    Thus $w$ is a weakly divergence free `vortex'. (c.f. Section 7 of \cite{HSSS}). Periodise $w$ by defining $\Lambda = \{(y_1, y_2) \in \mathbb{Z}^2 : y_1 + y_2 \text{ is even}\}$ and setting 
    \[
    u(x) = \sum_{y \in \Lambda} w(x-y)\, .
    \]
    Note that that $w$ is supported in one square of $\mathcal{S}^2$ and thus the periodization consists of filling half the squares of $\mathcal{S}^2$ with copies of $w$, while leaving the field identically equal to $0$ in the remaining squares. The ``filled'' and ``empty'' squares form likewise a chessboard pattern. 
    
    Even though $u$ is irregular, it has locally bounded variation and it is piecewise linear. There is thus a unique solution $\rho$ of (TE) with vector field $u$ and similarly the flux $\Phi$ of $u$ is well-defined. Its relevant property is that
    \begin{itemize}
    \item[(O)] The map $\Phi (t, \cdot)$ is Lipschitz on each square $S$ of $\mathcal{S}^2$ and $\Phi (\frac{1}{2}, \cdot)$ is a clockwise rotation of $90$ degrees of the ``filled'' $S$, while it is the identity on the ``empty ones''. In particular for every $j\geq 1$ $\Phi (\frac{1}{2}, \cdot)$ maps an element of $\mathcal{S}^1_j$ rigidly onto another element of $\mathcal{S}^1_j$. For $j=1$ we can be more specific. Each $S\in \mathcal{S}^2$ is formed precisely by $4$ squares of $\mathcal{S}^1_1$: in the case of ``filled'' $S$ the $4$ squares are permuted in a $4$-cicle clockwise, while in the case of ``empty'' $S$ the $4$ squares are kept fixed.  
    \end{itemize}
 Using this very last property it is therefore easy to see that
 \begin{equation}\label{e:refining}
 \rho(1/2, x) = 1 - \rho_{in}(2x)\, .
 \end{equation}
Likewise it is simple to use (O) to prove
\begin{itemize}
\item[(R)] If $\rho$ solves the transport equation \eqref{e:TE} with an initial data $\rho_{in}$ and $j\geq 1$, $\alpha$ are such that $\rho_{in}$ has average $\alpha$ on every $S\in\mathcal{S}^1_j$ with $j\geq 1$, then $\rho (\frac{1}{2}, \cdot)$ has also average $\alpha$ on $S\in \mathcal{S}^1_j$ 
\end{itemize}
 We define $b$ on $\mathbb R^2\times [0,2]$ in the following fashion. First of all $b(t, x) = u(x)$ for $0<t<1/2$ and $b(t, x) = u(2^nx)$ for $1-1/2^n<t<1-1/2^{n+1}$. For $1<t<2$, we let $b(t,x) = -b(2-t,x)$. 
    Note that \eqref{e:TE} has a unique solution $\rho$ on $[0, 1-1/2^n]$ because $b$ is a function of bounded variation. In particular this yields a unique solution on $[0,1]$. Moreover, using recursively the appropriately scaled version of \eqref{e:refining} we can readily check that $\rho (1-1/2^{2k}, x) = \rho_{in} (2^{2k} x)$ and $\rho (1-1/2^{2k+1},x) = 1 -\rho_{in} (2^{2k+1} x)$. In particular $\rho (t, \cdot) \rightharpoonup \frac{1}{2}$ as $t\to 1$. 
We can thus continue $\rho$ for $t\in [1,2]$ in two fashions, namely we set
\[
\rho'(t,x) = 
        \begin{cases} 
        \rho(t,x)\text{, for }0<t<1\\
        1/2\text{, for }1<t<2
        \end{cases}
\]
\[
\tilde \rho(t,x) = 
        \begin{cases} 
        \rho(t,x)\text{, for }0<t<1\\
        \rho(2-t,x)\text{, for }1<t<2
        \end{cases}\, .
\]
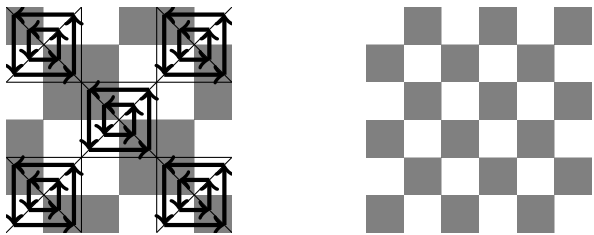
\begin{figure}[h]
\centering
\subfloat{
\begin{tikzpicture}
\clip(-1.5,-1.5) rectangle (1.5,1.5);
\foreach \x in {-1,...,2} \foreach \y in {-1,...,2}
    {
        \pgfmathparse{mod(\x+\y,2) ? "black" : "white"}
        \edef\nu{\pgfmathresult}
        \path[fill=\nu, opacity=0.5] (\x-1,\y-1) rectangle ++ (1,1);
    }
\draw (-1.5,0.5) -- (1.5,0.5);
\draw (-1.5,-0.5) -- (1.5,-0.5);
\draw (0.5, -1.5) -- (0.5, 1.5);
\draw (-0.5, -1.5) -- (-0.5, 1.5);
\foreach \a in {1,...,2}
    {
        \draw[ultra thick, ->] (\a/5,-\a/5) -- (\a/5,\a/5); 
        \draw[ultra thick, ->] (\a/5,\a/5) -- (-\a/5,\a/5);
        \draw[ultra thick, ->] (-\a/5,\a/5) -- (-\a/5,-\a/5);
        \draw[ultra thick, ->] (-\a/5,-\a/5) -- (\a/5,-\a/5);
    }
\foreach \a in {1,...,2}
    {
        \draw[ultra thick, ->] (\a/5+1,-\a/5+1) -- (\a/5+1,\a/5+1); 
        \draw[ultra thick, ->] (\a/5+1,\a/5+1) -- (-\a/5+1,\a/5+1);
        \draw[ultra thick, ->] (-\a/5+1,\a/5+1) -- (-\a/5+1,-\a/5+1);
        \draw[ultra thick, ->] (-\a/5+1,-\a/5+1) -- (\a/5+1,-\a/5+1);
    }
\foreach \a in {1,...,2}
    {
        \draw[ultra thick, ->] (\a/5+1,-\a/5-1) -- (\a/5+1,\a/5-1); 
        \draw[ultra thick, ->] (\a/5+1,\a/5-1) -- (-\a/5+1,\a/5-1);
        \draw[ultra thick, ->] (-\a/5+1,\a/5-1) -- (-\a/5+1,-\a/5-1);
        \draw[ultra thick, ->] (-\a/5+1,-\a/5-1) -- (\a/5+1,-\a/5-1);
    }
\foreach \a in {1,...,2}
    {
        \draw[ultra thick, ->] (\a/5-1,-\a/5-1) -- (\a/5-1,\a/5-1); 
        \draw[ultra thick, ->] (\a/5-1,\a/5-1) -- (-\a/5-1,\a/5-1);
        \draw[ultra thick, ->] (-\a/5-1,\a/5-1) -- (-\a/5-1,-\a/5-1);
        \draw[ultra thick, ->] (-\a/5-1,-\a/5-1) -- (\a/5-1,-\a/5-1);
    }
\foreach \a in {1,...,2}
    {
        \draw[ultra thick, ->] (\a/5-1,-\a/5+1) -- (\a/5-1,\a/5+1); 
        \draw[ultra thick, ->] (\a/5-1,\a/5+1) -- (-\a/5-1,\a/5+1);
        \draw[ultra thick, ->] (-\a/5-1,\a/5+1) -- (-\a/5-1,-\a/5+1);
        \draw[ultra thick, ->] (-\a/5-1,-\a/5+1) -- (\a/5-1,-\a/5+1);
    }
        \draw (-2,-2) -- (2,2);
        \draw (-2,2) -- (2,-2);
        \draw (-1.5,0.5) -- (-0.5, 1.5);
        \draw (1.5,0.5) -- (0.5,1.5);
        \draw (1.5,-0.5) -- (0.5, -1.5);
        \draw (-1.5,-0.5) -- (-0.5,-1.5);
\end{tikzpicture}%
}
\qquad
\qquad
\subfloat
{
\begin{tikzpicture}[scale=0.5]
\foreach \x in {-2,...,3} \foreach \y in {-2,...,3}
    {
        \pgfmathparse{mod(\x+\y,2) ? "white" : "black"}
        \edef\nu{\pgfmathresult}
        \path[fill=\nu, opacity=0.5] (\x-1,\y-1) rectangle ++ (1,1);
    }
\end{tikzpicture}%
}
\caption{Action of the flow from $t=0$ to $t=1/2$. The shaded region denotes the set $\{\rho=1\}$}
\end{figure}

{This is because we can `glue' weak solutions of \eqref{e:TE} to get another weak solution. More precisely, if $\rho_1$ and $\rho_2$ are weak solutions of \eqref{e:TE} defined for $0<t<1$ and $1<t<2$ respectively, and if both $\rho_1(t, \cdot)$ and $\rho_2(t, \cdot)$ weakly tend to the same limit as $t \to 1$, then 
\[
\rho(t,x) := 
        \begin{cases} 
        \rho_1(t,x)\text{, for }0<t<1\\
        \rho_2(t,x)\text{, for }1<t<2
        \end{cases}
\]
is also a weak solution of \eqref{e:TE}. Indeed, as $b$ is weakly divergence-free, for any $\phi \in C^\infty_c(\mathbb{R} \times \mathbb{R}^d)$ and $0<\beta<1$ we have
\begin{align}
    &\left|\int_0^\infty \int_{\mathbb{R}^d} \rho(t,x)(\partial_t \phi + (b \cdot \nabla) \phi) dx dt + \int_{\mathbb{R}^d} \rho_{in}(x)\phi(0,x) dx\right|
    \\&\leq \left|\int_0^{1-\beta} \int_{\mathbb{R}^d} \rho_1(t,x)(\partial_t \phi + (b \cdot \nabla) \phi) dx dt\right.\\
&\qquad\quad + \left.\int_{1+\beta}^\infty \int_{\mathbb{R}^d} \rho_2(t,x)(\partial_t \phi + (b \cdot \nabla) \phi) dx dt + \int_{\mathbb{R}^d} \rho_{in}(x)\phi(0,x) dx\right|
    \\&\qquad+ 2\beta\|D\phi\|_{L^\infty}\|b\|_{L^\infty}\|\rho\|_{L^1}
    \\& \leq \left|-\int_{\mathbb{R}^d} \rho_{in}(x)\phi(0,x) dx + \int_{\mathbb{R}^d} \rho_1(1-\beta,x) \phi(1-\beta,x) dx\right.\\
&\qquad\quad \left.- \int_{\mathbb{R}^d} \rho_2(1+\beta,x)\phi(1+\beta,x) dx + \int_{\mathbb{R}^d} \rho_{in}(x)\phi(0,x) dx\right|+ 2\beta\|D\phi\|_{L^\infty}\|b\|_{L^\infty}\|\rho\|_{L^1}
    \\& \leq \left|\int_{\mathbb{R}^d} \rho_1(1-\beta,x) \phi(1-\beta,x) dx - \int_{\mathbb{R}^d} \rho_2(1+\beta,x)\phi(1-\beta,x) dx \right|\\
&\quad+ 2\beta\|\partial_t \phi\|_{L^\infty} \|\rho_2(1+\beta,\cdot)\|_{L^1} + 2\beta\|D\phi\|_{L^\infty}\|b\|_{L^\infty}\|\rho\|_{L^1}
\end{align}
Thus, as $\beta \to 0$ and since both $\rho_1(t, \cdot)$ and $\rho_2(t, \cdot)$ weakly tend to the same limit as $t \to 1$, we get that $\rho$ is a weak solution of \eqref{e:TE}. 
}
\bigskip

\subsection{Step 2. Truncations.} We next construct two sequences of $BV$ vector fields converging to $b$. Both are simple and they are given by
\[
b_i^1(t,x) = 
    \begin{cases}
    b(t,x)\text{, for }0<t<1-1/2^{2i},\\
    0\text{, for }1-1/2^{2i}<t<1+1/2^{2i-2},\\
    b(t,x)\text{, for }1+1/2^{2i-2}<t<2
    \end{cases}
\]
\[
b^2_i(t,x) = 
    \begin{cases}
    b(t,x)\text{, for } 0<t<1-1/2^{2i},\\
    0 \text{, for }1-1/2^{2i}<t<1+1/2^{2i},\\
    b(t,x)\text{, for } 1+1/2^{2i}<t<2
    \end{cases}
\]
Let now $\rho^1_i$ and $\rho^2_i$ be the corresponding unique weak solutions of \eqref{e:TE} with initial data $\rho_{in}$. By construction both $\rho^1$ and $\rho^2$ coincide with $\rho = \rho' = \tilde{\rho}$ on the time interval $[0, 1-1/2^{2i}]$. Moreover for both we have
\begin{align*}
\rho^1_i (1+2^{-2i}, x) &=\rho (1-2^{-2i}, x) = \rho_{in} (2^{2i} x)\\
\rho^2_i (1+2^{-2i}, x) &= \rho (1-2^{-2i}, x) = \rho_{in} (2^{2i} x)\, .
\end{align*} 
Now, ${b}_i^2 (t,x) = b (t,x)$ for $t\geq 1+2^{-2i}$. Since $\rho^2_i (1+2^{-2i}, x) =  \rho_{in} (2^{2i} x) = \tilde{\rho} (1+2^{-2i}, x)$, we conclude that $\rho^2_i (t,x) = \tilde{\rho} (t,x)$ for $t\geq 1 + 2^{-2i}$. In particular we infer
\begin{equation}\label{e:second_solution}
\rho^2_i \rightharpoonup^\star \tilde{\rho}\qquad \mbox{in $L^\infty$.}
\end{equation}
When we come to $\rho^1_i (t,x)$ observe that, since $b (t,x)=0$ for $t\in [1+2^{-2i}, 1+ 2^{-2i+2}]$, we actually have $\rho^1_i (1+2^{-2i+2},x) = \rho_{in} (2^{2i} x)$. We can now use the (appropriately rescaled version of) (R) and conclude that
\begin{equation}\label{e:average}
\frac{1}{|S|}\int \rho^1_i (1+2^{-j}, x)\, dx = \frac{1}{2}\qquad \forall S\in \mathcal{S}^1_{2i-1}\quad
\mbox{and}\quad \forall j\leq 2i-2\, .
\end{equation}
In particular we conclude the same for $j=0$.
It is easy to see $\rho^1_i (2, \cdot) \rightharpoonup^\star \frac{1}{2}$ because of \eqref{e:average}. Indeed having fixed $\varphi\in C_c (\mathbb R^2)$, we can write
\[
\int \varphi (x) \rho^1_i (2,x)\, dx = \sum_{S\in \mathcal{S}^1_{2i-1}: S\cap {\mathrm{spt}}\, (\varphi)\neq \emptyset} \int_S \varphi (x) \rho^1_i (2,x)\, dx
\]
and hence estimate
\begin{equation}
\begin{split}
&\left|\int \varphi (x)\rho^1_i (2,x)\, dx - \frac{1}{2} \int \varphi (x)\, dx\right|\\
= &\left|\sum_{S\in \mathcal{S}^1_{2i-1}: S\cap {\mathrm{spt}}\, (\varphi)\neq \emptyset} \int_{S} \varphi (x) \left(\rho^1_i(2,x) - \frac{1}{|S|}\int_{S} \rho^1_i(2,x)\right)\right|\\
= &\left|\sum_{S\in \mathcal{S}^1_{2i-1}: S\cap {\mathrm{spt}}\, (\varphi)\neq \emptyset} \int_{S} (\varphi (x) -\varphi (x_S)) \left(\rho^1_i(2,x) - \frac{1}{|S|}\int_{S} \rho^1_i(2,x)\right)\right|
\\ \leq & \sum_{S\in \mathcal{S}^1_{2i-1}: S\cap {\mathrm{spt}}\, (\varphi)\neq \emptyset} \|D \varphi\|_{C^0} 2^{-2i+1} 2^{-2(2i-1)} \lesssim_{\varphi} \|D \varphi\|_{C^0} 2^{-2i+1}\, ,
\end{split}
\end{equation}
where we have $x_S$ to be the center of square $S$. {Now since for any $\epsilon>0$ any $f \in L^1$ can be written as $\varphi + h$ with $\varphi \in C^1$ and $\|h\|_{L^1} < \epsilon$, we get that 
\begin{equation}
\begin{split}
    \left|\int f (x)\rho^1_i (2,x)\, dx - \frac{1}{2} \int f (x)\, dx\right| & = \left|\int \varphi (x)\rho^1_i (2,x)\, dx - \frac{1}{2} \int \varphi (x)\, dx\right|\\
&\quad + \left|\int h (x)\rho^1_i (2,x)\, dx - \frac{1}{2} \int h (x)\, dx\right|\\
    & \lesssim_{\varphi} \|D \varphi\|_{C^0} 2^{-2i+1} + \|\rho^1_i (2,x) - {\textstyle{\frac{1}{2}}}\|_{L^\infty}\|h\|_{L^1}\\
& \lesssim_{\varphi} 2^{-2i+1} + \epsilon
    \end{split}
\end{equation}
Thus, as $\epsilon$ was arbitrary, we see that $\rho^1_i (2, \cdot) \rightharpoonup^\star \frac{1}{2}$.
}
So, any weak$^\star$ limit of a convergent subsequence of $\rho^1_i$ converges to a (backward) solution of the transport equation which is identically equal to $\frac{1}{2}$ at time $2$. We can now use the backward uniqueness for the transport equation with vector field $b$ on intervals $[1+\sigma, 2]$ for $\sigma>0$ (such uniqueness is guaranteed by the fact that the vector field $b$ is BV on $[1+\sigma,2]\times \mathbb R^2$), to conclude that such weak$^\star$ limit is identically equal to $\frac{1}{2}$ on $[1+\sigma, 2]$. In particular we conclude that $\rho^1_i$ converges weakly$^\star$ to $\rho'$.   

\subsection{Step 3. Regularization.} We now extend the vector fields $b^2_i$ and $b^1_i$ to times $t\not \in [0,2]$ by setting them identically $0$. We fix $i$ and a space-time compactly supported convolution kernel $\varphi$ and regularize both $b^1_i$ and $b^2_i$ to $b^1_{i,j}$ and $b^2_{i,j}$ setting $b^k_{i,j} := b^k_i * \varphi_{2^{-j}}$. Similarly, we regularize $\rho_{in}$ to $\rho_{j,in}:= \rho_{in} * \psi_{2^{-j}}$ for some space compactly supported convolution kernel $j$. Since each vector field $b^k_i$ belongs to $L^\infty ([0,2], BV\cap L^\infty (\mathbb  R^2))$, we can use Ambrosio's extension of the DiPerna-Lions theory to conclude that, for each fixed $i$ and $k$, the corresponding solutions $\rho^k_{i,j}$ of the transport equations with vector fields $b^k_{i,j}$ and initial data $\rho_{j,in}$ converge strongly in $L^1_{loc}$ to $\rho^k_i$. In particular we can select $j(i)$ such that
\[
\sum_{k=1}^2 \|\rho^k_{i,j(i)}-\rho^k_i\|_{L^1 ([0,2]\times [-2^i, 2^i]^2)} \leq 2^{-i}\, .
\]
We then set $b'_i = b^1_{i,j(i)}$, $\tilde{b}_i = b^2_{i,j(i)}$, $\rho'_i = \rho^1_{i,j(i)}$ and $\tilde{\rho}_i = \rho^2_{i,j(i)}$. Clearly, $\tilde{\rho}_i \rightharpoonup \tilde{\rho}$ and $\rho'_i \rightharpoonup \rho'$ in $L^1 ([0,T]\times U)$ for every bounded open set $U$.

\section{Proofs of Theorem \ref{t:main} and Theorem \ref{t:main-2}}\label{s:minor}

\subsection{Step 4. Compact supports.} Thus far the vector fields, the initial data and the solutions do not have compact support. However, in order to make them have compact supports we just proceed as follows. We modify $\rho_{in}$ to $\rho_{in} \mathbf{1}_{Q_N}$, where $Q_N = [-N,N]^2$ for some large natural number $N$. In order to truncate appropriately $b$ we need to act more carefully. For $t\in [1-2^{-i}, 1-2^{-i-1}]$ we substitute $b(t,x)$ with $b(t,x) \mathbf{1}_{Q_{N+2^{-i-1}}} (x)$. Observe that the choice of the sidelength of the square is made so to guarantee that the vector field remains divergence-free. We then keep the symmetric structure $b (t,x) = b(2-t,x)$ for the truncated field and we follow the same procedures of the previous steps. Note that the new regularized fields coincide with the old (nontruncated) ones in, say, $[0,2]\times Q_{N/2}$ and the initial data coincide with the old (nontruncated) ones in $Q_{N/2}$. Moreover the $L^\infty$ norm of all the fields is bounded uniformly by an absolute constant independent of $N$. In particular, for $N$ sufficiently large, the solutions of the transport equations with the truncated fields with truncated initial data coincide with the ones for the nontruncated fields and nontruncated initial data. We thus infer the same conclusions.

\subsection{Step 5. $W^{s,p}$ estimates.} We now show that our vector field $b$ of the previous section is in $W^{s, 1}_{loc}([0,T]\times \mathbb R^2)$ for every $s<1$. We'll make all our estimates on $B := [-1/2, 1/2]^2$ and $\Omega := [0,2] \times B$. Recall $b(t, x) = \pm u(2^{i+1}x)$ on $\mathcal{I}_i := (1-2^{-i}, 1-2^{-(i+1)}) \cup (1+2^{-(i+1)}, 1+2^{-i})$ and is identically $0$ elsewhere. Thus as,
\begin{equation} \begin{split}
    \|u(2^{i+1}\cdot)\|_{BV(B)} \lesssim 2^{di} \|w(2^{i+1}\cdot)\|_{BV(B)} & = 2^{di} (\|w(2^{i+1}\cdot)\|_{L^1(B)} + \|Dw(2^{i+1}\cdot)\|_{TV(B)}) \\ & = 2^{di} (2^{-di} \|w\|_{L^1(B)} + 2^{-(d-1)i} \|Dw\|_{TV(B)}) \\ & \leq 2^{di} 2^{-(d-1)i} \|w\|_{BV(B)} = 2^i \|w\|_{BV(B)}
    \end{split}
\end{equation}
The first inequality follows because there are approximately $2^{di}$ `little' vortices in $B$. Now,
\begin{equation}
    \|b \chi_{\mathcal{I}_i}\|_{BV(\Omega)} \leq C + \int_0^2 \|u(2^{i+1}x)\|_{BV(B)} dt \lesssim 1
\end{equation}
The constant $C$ comes form the `horizontal' jump part of the measure at $\partial \mathcal{I}_i$. Note this constant is indeed independent of $i$. By Gagliardo-Nirenberg we get for $0<s<1$,
\begin{equation}
    \|b \chi_{\mathcal{I}_i}\|_{W^{s, 1}(\Omega)} \lesssim \|b\chi_{\mathcal{I}_i}\|_{L^1(\Omega)}^{1-s} \|b\chi_{\mathcal{I}_i}\|_{BV(\Omega)}^{s} 
    \lesssim 2^{-i(1-s)}
\end{equation}
Thus,
\begin{equation}
    \|b\|_{W^{s, 1}(\Omega)} = \left\|\sum_{i=1}^\infty b \chi_{\mathcal{I}_i}\right\|_{W^{s, 1}(\Omega)} \lesssim \sum_{i=1}^\infty 2^{-i(1-s)} < +\infty\, .
\end{equation}
We leave to the reader the obvious modifications to deal with the truncation of $b$.

\subsection{Step 6. Making the field autonomous.}
    For any $a(t, x) \in L^\infty(\mathbb{R}_t \times \mathbb{R}^2; \mathbb{R}^2)$, we can define
    \begin{equation}
        f(a)(y) = 
            (1, a(y_0, y_1, y_2)) 
        \text{,} \ f(a) \in L^\infty(\mathbb{R}^{3}; \mathbb{R}^3)
    \end{equation}
    If we apply this transformation to the nonautonomous field defined in Step 1 of the previous section, we reach an autonomous field $v = f(b)$ and an initial density 
    \begin{equation}
        \theta_{in}(y_0, y_1, y_2) = \begin{cases}
        \rho_{in}(y_1, y_2) \text{ , for } -1 \leq y_0 \leq 0 \\
        0 \text{ , otherwise}
    \end{cases}
    \end{equation}
    Note that we again get 2 solutions for \eqref{e:TE}, namely:
    \begin{itemize}
        \item[i.] $\theta'(t,y)$ such that $\theta'(t,y) = 1/2$ for $2<t<3$ and $2<y_0<t$
        \item[ii.] $\tilde \theta(t,y)$ such that $\tilde \theta(t,y) = \rho_{in}(y_1, y_2)$ for $2<t<3$ and $2<y_0<t$
    \end{itemize}
We then apply the same procedure to all the nonautonomous fields constructed in the previous section to get an example which satisfies the requirements of Theorem \ref{t:main} with $d=2$. The extension to higher dimension is simple: in $\mathbb R^d = \mathbb R^3\times \mathbb R^{d-3}$ we just set the the components $v_j$ with $j\geq 2$ identically equal to $0$, while the remaining three components are made constant in the directions $y_3, \ldots, y_d$. Similarly the initial data is assumed constant along the directions $y_3, \ldots , y_d$. This then gives a noncompactly supported example: to pass to a compactly supported proceed as in the previous section.

The estimates obtained in Step 5 imply that $v\in W^{s,1}$ for every $s<1$. Fix now $s<1$ and select $\sigma\in (s,1)$. By interpolation we have
\[
\|v\|_{W^{s,p (s,\sigma)}} \leq C \|v\|_{L^\infty}^{1-\frac{s}{\sigma}}\|v\|_{W^{\sigma,1}}^{\frac{s}{\sigma}}
\]
for $\frac{1}{p (s,\sigma)} = \frac{s}{\sigma}$. Since we can take $\sigma$ arbitrarily close to $1$ we conclude that $v\in W^{s,p}$ for every $p<\frac{1}{s}$.

Observe next that, the vector fields $b^1_i$ and $b^2_i$ enjoy similar estimates, uniformly in $i$. Since the $\tilde{v}_i$ and  $v'_i$ are obtained from the latter through convolution with standard kernels and an application of $f$, the same uniform estimates are inherited by them. 

\bibliographystyle{plain}

\end{document}